\documentclass[12pt,a4paper]{article}
\usepackage{amsthm}
\usepackage{setspace}
\usepackage{plain}
\usepackage{xypic}
\usepackage{graphicx}
\usepackage{caption}
\begin{document}
\def\color#{}
\theoremstyle{plain}
	\newtheorem{theorem}{Theorem}
	\newtheorem{hypo}{Hypothesis}
	\newtheorem{lemma}{Lemma}
	\newtheorem{proposition}{Proposition}
	\newtheorem{corollary}{Corollary}
\def\N{{\bf N}}
\def\Z{{\bf Z}}
\def\R{{\bf R}}
\def\J{{\cal J}}
\def\K{{\cal K}}
\def\I{{\cal I}}
\def\words#1{\quad\hbox{#1}\quad}
\def\wwords#1{\qquad\hbox{#1}\qquad}
\def\<{\langle}\def\>{\rangle}
\def\dimf{\dim_{\rm b}}
{
    \catcode`\@=11
    \gdef\curl{\mathop{\operator@font curl}\nolimits}
}
\title{A Computational Study of a Data Assimilation Algorithm
for the Two-dimensional Navier--Stokes Equations}
\date{May 3, 2015}
\author{
Masakazu Gesho\thanks{
		Department of Mathematics,
		University of Wyoming,
		1000 E. University Ave, Dept. 3036,
		Laramie, WY 82071, USA.
		{\it email:}{\tt\ mgesho@uwyo.edu}} \and
Eric Olson\thanks{
		Department of Mathematics and Statistics,
		University of Nevada,
		Reno, NV 89557, USA. {\it email:}{\tt\ ejolson@unr.edu}} \and
Edriss S. Titi\thanks{
		Department of Mathematics, Texas A\&M University,
			3368--TAMU, College Station, TX 77843, USA.
		{\it email:}{\tt\ etiti@math.tamu.edu}}
	\thanks{
		The Department of Computer Science and Applied Mathematics,
		Weizmann Institute of Science,
		Rehovot 76100, Israel.
		{\it email:}{\tt\ edriss.titi@weizmann.ac.il}}
}

\maketitle

\begin{abstract}
We study the numerical performance of a continuous data assimilation (downscaling)
algorithm, based on ideas from feedback control theory, in the
context of the two-dimensional
incompressible Navier--Stokes equations.
Our model problem is to
recover an unknown reference solution, asymptotically in time, by using
continuous-in-time coarse-mesh
nodal-point observational measurements of the velocity field
of this reference solution (subsampling),
as might be measured by an array of weather vane anemometers.
Our calculations show that the required nodal observation density
is remarkably less that what is suggested by the analytical study; and is in fact comparable to the {\it number of numerically determining Fourier modes}, which was reported in an earlier computational study by the authors.
Thus, this method is computationally
efficient and performs far better than the analytical estimates
suggest.


\end{abstract}

\section{Introduction}

The goal of data assimilation is to provide a more accurate
representation of the current state of a dynamical system by
combining observational data with model dynamics.
This allows the influences of new data to be incorporated into
a numeric computation over time.
Data assimilation is widely used in the climate sciences, including
weather forecasting, environmental forecasting and hydrological forecasting.
Additional information and historical background may be found in
Kalnay \cite{kalnay2003} and references therein.

In 1969 Charney, Halem and Jastrow \cite{charney1969} proposed a
method of continuous data assimilation in which observational
measurements are directly inserted into the mathematical model
as it is being integrated in time.
To fix ideas, let us suppose that
the evolution of $u$ is governed by the dynamical system
\begin{equation}
\label{refsol}
	{du\over dt} = {\cal F}(u),\qquad u(t_0)=u_0
\end{equation}
and the observations of $u$ are given by the time series
$p(t)=P u(t)$ for $t\in[t_0,t_*]$, where $P$ is an orthogonal
projection onto the low modes.
In this context, the method proposed in \cite{charney1969} for
approximating $u$ from the observational data is to solve
for the high modes
\begin{equation}
\label{alg1}
{d q\over dt} = (I-P) {\cal F}(q+p),\qquad q(t_0)=q_0
\end{equation}
where $q_0$ is an arbitrarily chosen initial condition
and $q+p$ represents the resulting approximation of $u$.
Note that if $q_0=(I-P)u_0$ then $p+q=u$ for all time;
however, data assimilation is applied when $u_0$ is not known.

Algorithm (\ref{alg1}) was studied by Olson and Titi in
\cite{olson2003} and \cite{olson2008}
for the two-dimensional
incompressible Navier--Stokes equations
\begin{plain}\begin{equation}\label{nseq}\left\{\eqalign{
	{\partial u\over \partial t} -\nu \Delta u
	+(u\cdot \nabla)u+\nabla p = f\cr
	\nabla\cdot u = 0
	}\right.
\end{equation}\end{plain}%
on the
domain $\Omega=[0,L]^2$, equipped with periodic
boundary conditions and zero spatial average
with initial condition $u(x,t_0)=u_0(x)$ for $x\in [0,L]^2$.
Observational measurements were represented by
$P=P_h$, where $P_h$ is the orthogonal projection onto
the Fourier modes $\exp(2\pi  ik\cdot x /L)$ with wave
numbers $k\in\Z^2\setminus\{0\}$ such that $0<|k|\le L/h$.
Here $\nu>0$ is the kinematic viscosity,
$p(x,t)$ is the pressure and $f(x)$ is a time-independent
body force with zero spatial average acting on the fluid.
For simplicity, it was
assumed, as we shall here, that $\nabla\cdot f=0$.

The two-dimensional incompressible Navier--Stokes equations are
amen\-able to mathematical analysis while
at the same time they posses non-linear dynamics
similar to the partial differential equations that govern
realistic physical phenomenon.
Using the functional notation of
Constantin and Foias \cite{constantin1988}, see also
Temam~\cite{temam1983} or Robinson \cite{robinson2001},
write (\ref{nseq}) in the form (\ref{refsol}) by setting
\begin{equation}
\label{NSEQ}
    {\cal F}(u)=-\nu A u - B(u,u)+f
\end{equation}
where $A$ and $B$
are the continuous extensions of
the operators given by
$$A=-P_\sigma \Delta u\wwords{and}
B(u,v)=P_\sigma (u\cdot \nabla v)$$
when $u,v$ are smooth divergence-free $L$-periodic functions,
and $P_\sigma$ is the Leray--Helmholtz projector.
We recall that
$$
	P_\sigma(u)= \sum_{k\in\Z^2\setminus\{0\}}
		\bigg\{ u_k - {k\cdot u_k\over |k|^2}k\bigg\}
		\exp(2\pi i k\cdot x/L)
$$
and also that $A\colon V^1\to V^{-1}$ and $B\colon V^1\times V^1\to V^{-1}$
where
$V^\alpha$ is the closure of~${\cal V}$, the space of zero-average
$\R^2$-valued
divergence-free $L$-periodic trigonometric polynomials, with respect
to the norm
$$
	\|u\|_{V^{\alpha}}
		= L^2 \sum_{k\in\Z^2\setminus\{0\}} |k|^{2\alpha} |\hat u_k|^2.
$$
For notational convenience we shall write $V=V^1$ throughout
the remainder of this paper.

Consider
the data assimilation method given by (\ref{alg1}).
Using the theory of determining modes it was shown as
Theorem 1.5 in \cite{olson2003} that
if $h$ satisfies
\begin{equation}\label{alg1thm}
    {L^2\over h^2}\ge c_1 G
\end{equation}
then $\big\|u(t)-p(t)-q(t)\big\|_V \to 0$,
exponentially fast as $t\to\infty$.
Here $c_1$ is a universal constant and
$G=({L/2\pi \nu})^2 \|f\|_{L^2}$ is the Grashof number.
Computations in \cite{olson2008} considered
a fixed forcing function $f=f_{121}$ scaled to obtain
different values for $G$.
In that work the subscript $121$ was used to indicate
that the force was supported on an annulus around $k^2=121$
in Fourier space.
More details on $f$ are provided by equation (\ref{force})
and Figure~\ref{frcinit} below.
For Grashof numbers between $500\,000$ and $60\,000\,000$
it was shown that the projection $P_h$ onto the lowest
$80$ Fourier modes was necessary and sufficient
to ensure numerically that $\big\|u(t)-p(t)-q(t)\big\|_V\to 0$, as
$t\to\infty$.
Since the rank of $P_h$ scales as $L^2/h^2$, this is
significantly less than the millions suggested by the
analytical bound (\ref{alg1thm}).
Thus, the data assimilation algorithm given by~(\ref{alg1}) performs
far better than the analysis suggests.
Note, however, that this algorithm
is not suitable when the observations
are given by nodal measurements of the velocity field.

An approach to data assimilation, for dynamics governed by
equations (\ref{nseq}) or equivalently (\ref{NSEQ}), which is
applicable to
nodal observations,
was introduced and analyzed by Azouani, Olson and Titi in \cite{azouani2014}.
Let
$I_h$ be a general interpolant observable satisfying the
approximation identity inequality
\begin{equation}
\label{genpoin}
	\|u-I_h(u)\|_{L^2}^2\le \gamma_1 h^2 \| u\|_{H^1}^2
			+ \gamma_2 h^4 \|u\|_{H^2}^2.
\end{equation}
Given $I_h(u(t))$ for $t\in[t_0,t_*]$,
solve
\begin{equation}
\label{alg2}
	{dv\over dt} = {\cal F}(v)+\mu P_\sigma \big(I_h(u)-I_h(v)\big),
	\qquad v(t_0)=v_0,
\end{equation}
where $v_0$ is an arbitrary initial condition.
The constant $\mu$ is a relaxation (nudging) parameter which
controls the strength of the feedback control (nudging term).
In particular, the nudging term pushes the large spatial
scales of the approximating solution $v$ toward those of
the reference solution $u$ while the viscosity stabilizes
and dissipates the fine spatial scales and any spillover into
the fine scales caused by the nudging term.
It follows from Theorem~2 equation~(39) of~\cite{azouani2014}
that if $h$ and $\mu$ satisfy
\begin{equation}\label{alg2thm}
{L^2\over h^2}\ge { c_0 L^2\mu\over \nu}\ge c_2 G\big(1+ \log (1+G)\big),
\end{equation}
then
$\|u(t)-v(t)\|_V\to 0$ as $t\to\infty$.
Here $c_2$ is a universal constant and
$c_0$ is a constant depending only on $\gamma_1$ and $\gamma_2$
of (\ref{genpoin}).

The number of nodal measurements needed to uniquely
determine a
solution to the two-dimensional Navier--Stokes equations, as $t\to\infty$,
was found by Foias and Temam
in~\cite{foias1984} and further refined in
Jones and Titi~\cite{jones1993}.
Up to a logarithmic correction, the analytic bounds
given by (\ref{alg2thm}) are the same as those given in \cite{jones1993}.
In this paper we check the numerical performance of the data
assimilation algorithm (\ref{alg2}) using nodal measurements
given by $I_h$
for the same body forcing $f$
considered
in \cite{olson2008} scaled so that $G=2\,500\,000$.

Let $Q_i$ be disjoint squares that cover $[0,L]^2$
with centers
$x_i$ and sides of length $h=L/K$, where $K^2=N$.
An interpolant operator based on the nodal measurements $u(x_i,t)$,
for $i=1,2,\ldots,N$, and $t\in [t_0,t_*]$,
which satisfies (\ref{genpoin}) is
\begin{equation}
\label{nodes}
    I_h(u)(x,t)={\cal I}_h(u)(x,t)-{1\over L^2}
    \int_{[0,L]^2} {\cal I}_h(u)(x,t)\,dx,
\end{equation}
where
\begin{equation}
    {\cal I}_h(u)(x,t)=\sum_{i=1}^{N} u(x_{i},t)\chi_{Q_{i}}(x).
\end{equation}

We also consider the smoothed interpolation
\begin{equation}
\label{snodes}
	\tilde I_h(u)(x,t)=\widetilde{\cal I}_h(u)(x,t)-{1\over L^2}
	\int_{[0,L]^2} \widetilde{\cal I}_h(u)(x,t)\,dx,
\end{equation}
where
\begin{equation}
	\widetilde{\cal I}_h(u)(x,t)
		=\sum_{i=1}^{N} u(x_{i},t)(\rho_{\epsilon} *\chi_{Q_{i}})(x),
\end{equation}
and $\rho_\epsilon(x)=\epsilon^{-2}\rho(x/\epsilon)$ with
\begin{plain}\begin{equation}
\rho(\xi)=\cases{
	K_0\exp\bigg(\displaystyle{1\over 1-\xi_1^2}+{1\over 1-\xi_2^2}\bigg)
		&for $|\xi_1|,|\xi_2|<1$\cr
	0&otherwise,
	}
\end{equation}\end{plain}%
and
$$
	K_0^{-1}=
	\int_{-1}^1
	\int_{-1}^1
\exp\bigg(\displaystyle{1\over 1-\xi_1^2}+{1\over 1-\xi_2^2}\bigg)
		\,d\xi_2\, d\xi_1.
$$

To make the smoothing scale compatible with the resolution parameter
we take $\epsilon=\eta h$ for some fixed $\eta>0$.
An analysis of a smoothed interpolant similar to $\tilde I_h$
appears in Appendix A of \cite{azouani2014} for $\eta=0.1$ and
shows that $\tilde I_h$ satisfies~(\ref{genpoin}).  Note also
that $\tilde I_h\to I_h$ as $\eta\to 0$.
When $\eta$ is between $0$ and $1$ the convolution
reduces the high-frequencies that would otherwise be present
in the Fourier series representation of
the characteristic
function $\chi_{Q_i}$.  Values of $\eta$ greater than $1$
blur nearby nodal measurements together.
This further downscaling could be useful in the presence of
noisy measurements, see for example \cite{bessaih2015};
however, the measurements studied here will be
error free.
In this work we vary $\eta$ between $0.1$ and $2.0$ and find
that $\eta=0.7$ leads to
near optimal performance for our data assimilation
experiments.

In particular, our results show that $\|u(t)-v(t)\|_V\to 0$,
as $t\to\infty$, when the resolution $K$ of the observational
measurements satisfies $K\ge 8$ and $\mu$ and $\eta$ are
appropriately chosen.
Moreover, if $K\ge 9$ and $\eta=0.7$, there is a wide range of
values for $\mu$ such that the algorithm works well.
Since $64$ and $81$ nodes are comparable in resolution
to $80$ Fourier modes, the numerical efficiency of algorithm
(\ref{alg2}), using nodal measurements, is comparable
to algorithm (\ref{alg1}) using Fourier modes
\cite{olson2003,olson2008}.

Rigorous mathematical analysis of the method of data
assimilation studied computationally in this paper
has recently been generalized to B\'enard convection by
Farhat, Jolly and Titi~\cite{farhat2015} where it was shown that only observational measurements of the velocity field is sufficient to recover the full reference solution, i.e., the velocity field and the temperature. Inspired by~\cite{farhat2015} Farhat, Lunasin and Titi~\cite{FLT} have recently improved the algorithm studied here, i.e.~the one introduced in~\cite{azouani2014}, by showing that it is sufficient to use observational measurements of only one
component of the velocity field to recover the full reference solution. Further implementation of this algorithm,  for the subcritical
surface quasi-geostrophic equation, has recently been established by Jolly, Martinez and Titi~\cite{jolly2015}.
The algorithm studied here is also closely related
the 3DVAR data assimilation method developed by Bl\"omker,
Law, Stuart and Zygalakis~\cite{blomker2013} for the
Navier--Stokes equations and by Law, Shukla and
Stuart~\cite{law2014} for Lorenz equations.

This paper is organized as follows.
Section 2 describes the physical parameters, forcing and
initial conditions used to
generate the reference solution to the two-dimensional
incompressible Navier--Stokes equations that we will be
observing through nodal measurements of the velocity field.
Section 3 reports our computational results, and
section 4 gives details of our computational methods.
The last section concludes that
data assimilation of nodal measurements, by means of
equation~(\ref{alg2}), as studied in this paper works
computationally just as efficiently
as equation~(\ref{alg1}) with Fourier modes.

\section{The Reference Solution}
To focus on how the smoothing and resolution of the observational
measurements affect algorithm (\ref{alg2}), as
well as how to optimize the
value of the relaxation (nudging) parameter $\mu$,
we fix the viscosity and the
size of the periodic box so that
$$
\nu=0.0001
\wwords{and}
L=2\pi
$$
for the remainder of this paper.
We further perform all our simulations using the same reference solution
$u(t)$ to the two-dimensional incompressible Navier--Stokes equations.
As shown in \cite{olson2003} the spatial structure
of the function $f$ used to force the reference solution can
have a significant effect on data assimilation.
Therefore, to allow comparison with previous results, we
use the exact same forcing function
defined in \cite{olson2003}, and further
studied in \cite{olson2008}, scaled such that $G=2\,500\,000$
in our present computations.

This function $f$ is supported on the annulus in Fourier space
with wave numbers $k$ such that $110\le k^2\le 132$.
In particular,
\begin{equation}
\label{force}
f(x)=\sum_{110\le k^2\le 132} \hat f_k \exp(ik\cdot x)
\end{equation}
with
$\hat f_k=\overline{\hat f_{-k}}$
and
$k\cdot \hat f_k=0$,
where the values of $\hat f_k$ are given by
Table~1 in~\cite{olson2008} scaled to obtain
the desired Grashof number.
Note that this forcing is at length scales of about
$1/11$-th the size of the periodic box.
This fact is further reflected in the level curves of
$$\curl f=\curl(f_1,f_2)={\partial f_2\over \partial x}-{\partial f_1\over \partial y}$$
depicted
in Figure \ref{frcinit} on the left.

\def\mysize{3in}%
\begin{figure}[h!]
	\caption{Left are contours of $\curl f$; right
		are vorticity contours
		of $u_0$.}
	\label{frcinit}
	\vskip-1cm
	\centerline{
	\includegraphics[width=\mysize]{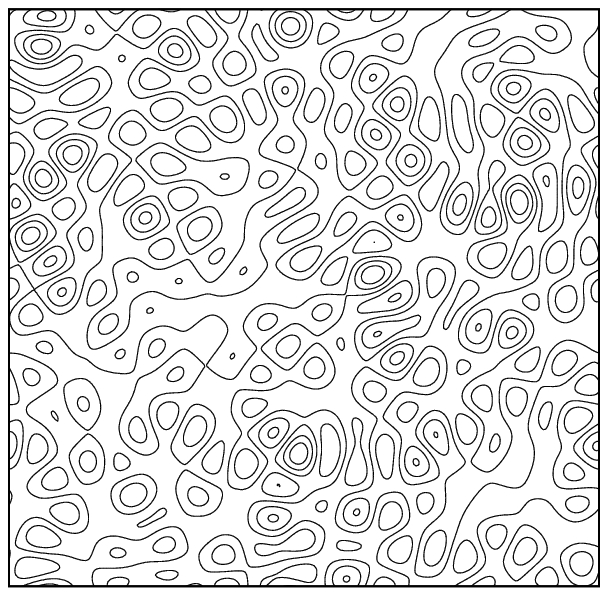}\hskip-1.2cm
	\includegraphics[width=\mysize]{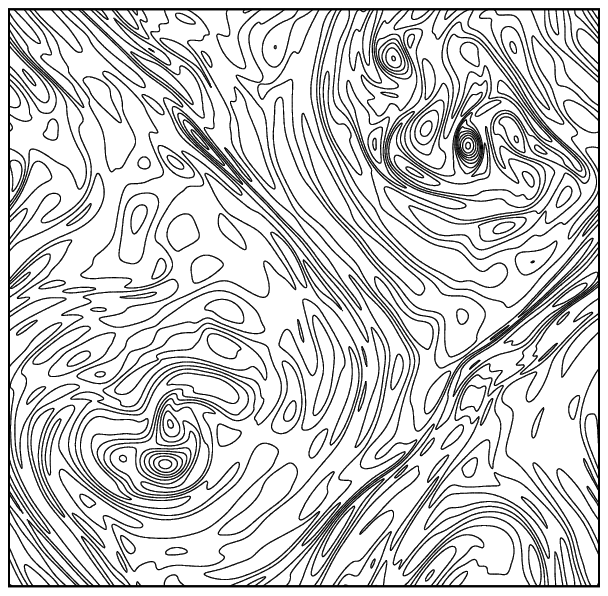}
}
\vskip -1cm
\end{figure}

The initial condition $u_0$ used for our data assimilation experiments
was obtained by solving (\ref{nseq}) with zero initial condition at time
$t=-25000$ until time $t=0$.
In terms of eddy turnover times, this ensures that more than
$500$ eddy turnovers have occurred before reaching $t=0$.
Integrating for this length of time ensures the
initial condition $u_0$ lies close to the global attractor
and therefore reflects the energetics of the forcing $f$.
In particular, the way in which we
initialized the solution at time $t=-25000$ is unimportant.

The vorticity contours $
	\omega_0=\curl u_0$ of the initial condition
$u_0$ are depicted in Figure \ref{frcinit} on the right.
While the forcing $f$ contains no Fourier modes with
wave numbers $k$ such that $|k|<10$, the initial condition
$u_0$ clearly possesses two large eddies the size of the box.
These large box-filling eddies apparently result from the
inverse cascade of energy in the two-dimensional
Navier--Stokes equations.  This can be seen more clearly
by examining the energy spectrum of the reference solution.

\begin{figure}[h!]
	\caption{Time averaged energy spectrum of the reference
	solution.}
	\label{spec}
	\centerline{
	\includegraphics[width=4in]{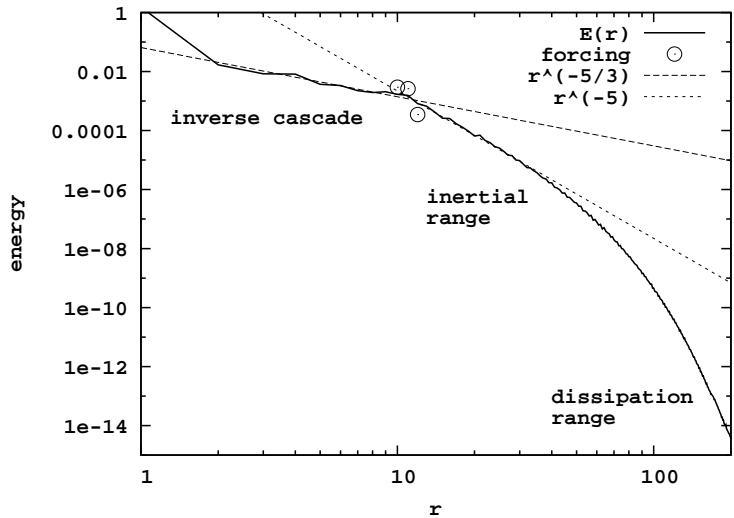}}
\end{figure}

Given a solution $u(t)$ to the two-dimensional
Navier--Stokes equations for $t\in[0,T]$ where $T=25\,000$
define the average energy spectrum as
$$
	E(r)={4\pi^2\over T} \int_{0}^T
		\sum_{ k\in\J_r} |\hat u_k(t)|^2 dt
\words{where}
	u(t)=\sum_{k\in\Z^2\setminus\{0\}} \hat u_k(t) e^{ik\cdot x}
$$
and $\J_r=\big\{\, k\in\Z^2 : r-0.5<|k|\le r+0.5\,\big\}$.
For the reference solution described above with initial
condition $u_0$ and time $t_0=0$, the average energy
spectrum appears in Figure \ref{spec}.
While we do not see the Kraichnan scaling of $k^{-3}$
in the inertial range, we do see the Kolmogorov $k^{-5/3}$
scaling in the inverse cascade.
Such spectra have been observed in other numerical
experiments, see, for example,
Xiao, Wan, Chen and Eyink \cite{xiao2009}.
In particular, the inverse cascade appears responsible for
the box-filling eddies observed in the initial condition $u_0$
which persist for all times $t>0$.

We compute the eddy turnover time for the reference
solution as
$$
	\tau=4\pi^2\sum_{r=1}^\infty r^{-1} E(r)
		/\Big(\sum_{r=1}^\infty E(r)\Big)^{3/2}
\approx 30.8
$$
and conclude that our averages have been computed over
$T/\tau\approx 812$ eddy turnovers.
The spectrum of $\tau f$, which also has units of
energy, has been plotted
in Figure \ref{spec} as three circles
to illustrate where the forcing lies in relation
to the energy spectrum.  Note that the forcing
exactly divides the energy spectrum between the part which
scales as $k^{-5/3}$ and the part which scales as $k^{-5}$.

Having, to some extent, described
the reference solution that will be used in our numerical
experiments,
we now turn to our main point of study, the data assimilation of
nodal measurements of the velocity field.

\section{Nodal Observations of Velocity}

We consider nodal observations $u(x_i,t)$, for $i=1,\ldots,N$,
of the reference solution $u$,
that was computed according to the incompressible two-dimensional
Navier--Stokes equations (\ref{nseq}) and initialized with
$u_0$ as described in Figure \ref{frcinit},
and interpolate these measurements
using the operator $I_h$ defined by (\ref{nodes}).
The resulting equations for
the approximating solution $v$ may be written as
\begin{equation}\label{alg2sys}
	{dv\over dt}
		+\nu A v + B(v,v) = f - \mu P_\sigma (I_h(u)-I_h(v))
\end{equation}
where $v$ is initialized as $v_0=0$.
Note that only the observations $I_h(u)$ of the reference
solution $u$ enter into the equations for computing for $v$.
Also note that
$\|u(0)-v(0)\|_V =\|u_0\|_V\approx 1.946.$
Our goal now is to choose the resolution parameter~$h$
and the relaxation (nudging) parameter~$\mu$
in such a way that $\|u(t)-v(t)\|_V\to 0$,
numerically, as $t\to\infty$.

\begin{figure}[h!]
	\caption{The error $\|u(t)-v(t)\|_V$ versus $t$ for $h=0.6981$.}
	\label{decaymu}
	\centering
	\includegraphics[width=4in]{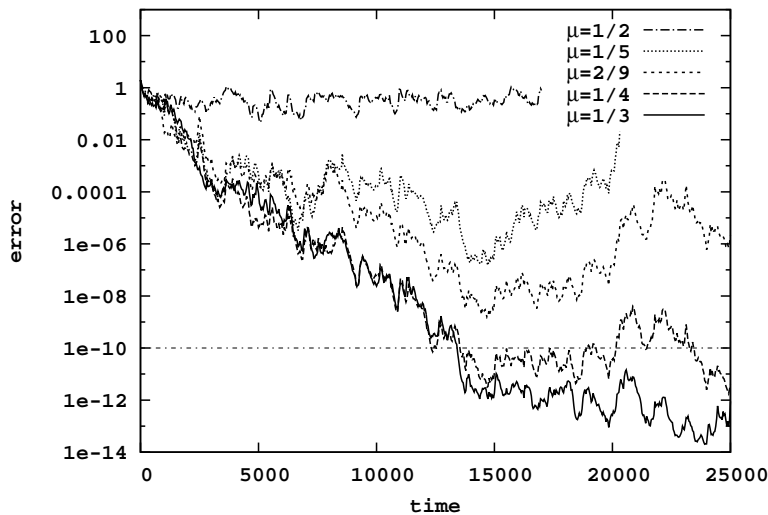}
\end{figure}

As discussed in \cite{azouani2014}, if $\mu$ is too small,
the feedback control (nudging term) will be too
weak to ensure the approximating solution converges to
the reference solution.
If $\mu$ is too large, then spill over into the fine scales
becomes significant and again prevents recovery of the
reference solution.  Figure \ref{decaymu} illustrates
each of these possibilities
for $h=L/K$, where $K=9$, using different values of $\mu$.

When $\mu=1/2$ the relaxation (nudging) parameter is too large
for the approximating solution to converge to the
reference solution, and when $\mu=1/5$ it is too small.
However, the intermediate value $\mu=1/3$ works
with the error represented by $\|u(t)-v(t)\|_V$
falling below $10^{-10}$ by $t=13\,417.8$.
Note that, since the double-precision floating-point numbers
used to represent the Fourier modes of $u$ and $v$ on the computer
have only 15 digits of precision, we don't expect convergence
of $\|u(t)-v(t)\|_V$ to exact zero over time.

For $\mu=1/4$ the error falls below $10^{-10}$
at $T=12\,327.1$,
however,
it rises again and it is not clear whether
after $T=23\,463.9$ the error
finally stays below $10^{-10}$
or not.
The value $\mu=2/9$ shows an even more irregular
pattern where $\|u(t)-v(t)\|_V$
exhibits a period of decay
followed by a period of growth that covers six orders of magnitude.
Fortunately, most of our parameter choices avoid these borderline
cases and the corresponding error either converges
towards zero and stays below $10^{-10}$ or shows few signs of
converging and stays well above $10^{-10}$.

To determine the values of $h$ and $\mu$ for which it is possible to
recover the reference solution to within numerical roundoff error
we fix $\epsilon>0$ and define
$$
	T_{\rm max}
		=\sup\big\{\, t\in [0,T]: \big\|v(t)-u(t)\big\|_V
			\ge\epsilon\,\big\}
$$
and
$$
	T_{\rm min}
		=\inf\big\{\, t\in [0,T]: \big\|v(t)-u(t)\big\|_V
			\le\epsilon\,\big\}
$$
Let
$T_{\rm max}=0$,
when the supremum is over the empty set, and
$T_{\rm min}=\infty$
when the infimum is over the empty set.
When $\|v(T)-u(T)\|_V\ge\epsilon$ further
set $T_{\rm max}=\infty$ to ensure
$T_{\rm max}\ge T_{\rm min}$.
Inspired by Figure \ref{decaymu} we also define
$$\varepsilon_{\rm avg} = {1\over T-T_0}\int_{T_0}^T \|v(t)-u(t)\|_V dt
$$
and take
$\epsilon=10^{-10}$, $T=25\,000$ and $T_0=2T/3$
for our numerics.

\begin{table}[h!]
	\caption{Data assimilation using $I_h$.}
	\centering				
	\begin{tabular}{r | r r l | r r l} 	
	\hline\hline 			
	& \multicolumn{3}{c|}{$K=9$}
	& \multicolumn{3}{c}{$K=10$} \\
	$\mu$ & $T_{\rm min}$ & $T_{\rm max}$ & $\varepsilon_{\rm avg}$
	& $T_{\rm min}$
	& $T_{\rm max}$ & $\varepsilon_{\rm avg}$ \\ [0.5ex]
	\hline	
  0.0625&    $\infty$&    $\infty$&  $1.1$&  
             $\infty$&    $\infty$&  $8.8\times 10^{-1}$\\ 
   0.125&    $\infty$&    $\infty$&  $3.9\times 10^{-1}$&  
          13754.5& 17094.3&  $4.1\times 10^{-12}$\\ 
   0.154&    $\infty$&    $\infty$&  $2.4\times 10^{-1}$&  
           4331.1&  4432.6&  $5.4\times 10^{-14}$\\ 
   0.167&    $\infty$&    $\infty$&  $2.4\times 10^{-1}$&  
           3965.0&  4187.7&  $3.1\times 10^{-14}$\\ 
   0.182&    $\infty$&    $\infty$&  $1.5\times 10^{-1}$&  
           3320.1&  3320.1&  $2.8\times 10^{-14}$\\ 
     0.2&    $\infty$&    $\infty$&  $3.5\times 10^{-2}$&  
           2825.2&  2825.2&  $2.8\times 10^{-14}$\\ 
   0.222&    $\infty$&    $\infty$&  $1.9\times 10^{-5}$&  
           2870.1&  2870.1&  $2.5\times 10^{-14}$\\ 
    0.25& 12327.1& 23463.9&  $3.7\times 10^{-10}$&  
           2701.0&  2701.0&  $2.5\times 10^{-14}$\\ 
   0.286& 12275.6& 13466.2&  $1.3\times 10^{-12}$&  
           2581.4&  2598.6&  $2.1\times 10^{-14}$\\ 
   0.333& 13417.8& 13417.8&  $1.4\times 10^{-12}$&  
           2601.8&  2601.8&  $2.0\times 10^{-14}$\\ 
     0.4&    $\infty$&    $\infty$&  $1.3\times 10^{-1}$&  
           3008.4&  3008.4&  $2.1\times 10^{-14}$\\ 
     0.5&    $\infty$&    $\infty$&  $3.2\times 10^{-1}$&  
           4564.3&  4564.3&  $2.4\times 10^{-14}$\\ 
     0.6&    $\infty$&    $\infty$&  $6.7\times 10^{-1}$&  
           8604.5&  9964.6&  $5.7\times 10^{-14}$\\ 
     0.7&    $\infty$&    $\infty$&  $1.5$&  
           $\infty$&    $\infty$&  $5.3\times 10^{-2}$\\ 
	[1ex]
	\hline	
	\end{tabular}
	\label{hmaxu}
\end{table}

Table \ref{hmaxu} shows the results of our computational
experiments.
Runs with $K=8$ were also performed, however, no value of
$\mu$ yielded a finite value for $T_{\rm max}$ or $T_{\rm min}$
or even an approximation for which the error
$\|u(t)-v(t)\|_V$ fell below $10^{-2}$.
We conclude
that $K=9$ is the minimal resolution for which there exists
a $\mu$ such that the error tends toward zero.
At this minimal resolution
only a narrow range of values for $\mu$
near $1/3$, result in an error which falls below $10^{-10}$.
When $K=10$, there is a much
greater range of corresponding values for $\mu$
that work well.
In fact, when $K=10$ all values of $\mu$ between $1/6$ and $1/2$
led to corresponding approximations $v(t)$ such that
$\varepsilon_{\rm avg}\approx 3\times 10^{-14}$.
Since the double-precision floating-point numbers
used to represent the Fourier modes of $u$ and $v$ on the computer
have only 15 digits of precision, the fact that the error
can approach $10^{-14}$ is remarkable.
We again note, as is consistent with the analysis
in \cite{azouani2014}, our numerical experiments do
not recover the reference solution
if $\mu$ is too small or too large.

\begin{figure}[h!]
	\caption{$T_{\rm max}$
		versus $\mu$ for $\epsilon=10^{-10}$, $T=25000$ and $K=10$.}
	\label{muopt10}
	\centering
	\includegraphics[width=0.7\textwidth]{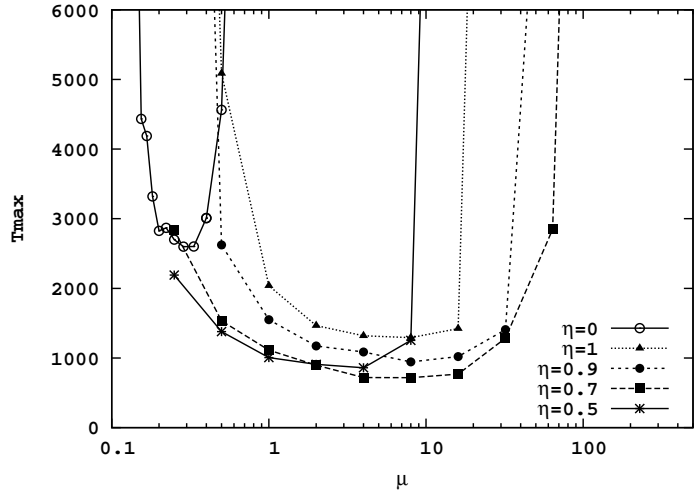}
\end{figure}

Next consider the family of smoothed interpolants $\tilde I_h$
for different values of $\eta$.
Figure \ref{muopt10} plots $T_{\rm max}$ versus $\mu$.
When $\eta$ is near $0.7$ we find that values of $\mu$
between $1/4$ and $64$ all lead to approximations
such that $\|u(t)-v(t)\|_V\le 10^{-10}$ for large enough $T$.
Thus, smoothing with $\eta$ near $0.7$ leads to a
significantly wider range of values for $\mu$ such that
the data assimilation algorithm can be used to
recover the reference solution.
Note when $\eta=2$ and $K=10$ that no values of $\mu$
led to the convergence of the approximating solution to
the reference solution over time.

Having found good values for $\eta$ we continue our
numerical study by fixing $\eta=0.7$ and
varying $\mu$ for different resolutions $h=L/K$, where $K=8$, 9 and 10.
As before $\epsilon=10^{-10}$, $T=25\,000$ and $T_0=2T/3$.
The computational results given in Table \ref{hmaxusmooth}
show that observational measurements with a resolution given
by $K=8$ can now lead to an approximate
solution which converges to the reference solution over time.
Moreover, the accuracy of the resulting approximations also
improve compared to the non-smoothed case.
Note that
we have omitted reporting
$T_{\rm min}$ in Table \ref{hmaxusmooth} since in all cases
$T_{\rm min}$ was equal or nearly equal to $T_{\rm max}$.

\begin{table}[h!]
    \caption{Data assimilation using $\tilde I_h$ where $\eta=0.7$.}
    \centering
    \begin{tabular}{r | r l |r l| r l}
    \hline\hline
    & \multicolumn{2}{c|}{$K=8$}
    & \multicolumn{2}{c|}{$K=9$}
    & \multicolumn{2}{c}{$K=10$} \\
    $\mu$ & $T_{\rm max}$ & $\varepsilon_{\rm avg}$
		& $T_{\rm max}$ & $\varepsilon_{\rm avg}$
    & $T_{\rm max}$ & $\varepsilon_{\rm avg}$ \\ [0.5ex]
    \hline
0.25&   $\infty$& $5.1\times 10^{-1}$ &
        $\infty$& $9.9\times 10^{-2}$ &
       2834.2& $2.4\times 10^{-14}$ \\
 0.5&   $\infty$& $6.5\times 10^{-2}$ &
       2570.1& $2.4\times 10^{-14}$ &
       1534.2& $1.7\times 10^{-14}$ \\
   1&  2817.8& $1.1\times 10^{-13}$ &
       1686.5& $1.8\times 10^{-14}$ &
       1112.7& $1.6\times 10^{-14}$ \\
   2&  2527.1& $2.6\times 10^{-14}$ &
       1232.6& $1.7\times 10^{-14}$ &
        897.5& $1.6\times 10^{-14}$ \\
   4&  2013.1& $2.1\times 10^{-14}$ &
       1092.4& $1.7\times 10^{-14}$ &
        718.9& $1.6\times 10^{-14}$ \\
   8&  2191.6& $2.1\times 10^{-14}$ &
       1124.8& $1.7\times 10^{-14}$ &
        717.8& $1.6\times 10^{-14}$ \\
  16&  4137.3& $4.7\times 10^{-14}$ &
       1360.4& $1.7\times 10^{-14}$ &
        769.7& $1.6\times 10^{-14}$ \\
  32&   $\infty$& $3.0\times 10^{-1}$ &
       2752.9& $2.1\times 10^{-14}$ &
       1284.2& $1.7\times 10^{-14}$ \\
  64&   $\infty$& 1.4 &
        $\infty$& $4.9\times 10^{-1}$ &
       2848.6& $2.1\times 10^{-14}$ \\
 128&   $\infty$& 2.6 &
        $\infty$& 1.9 &
        $\infty$& 1.1 \\
    [1ex]
    \hline
    \end{tabular}
    \label{hmaxusmooth}
\end{table}

\begin{figure}[h!]
	\caption{$T_{\rm max}$
		versus $\mu$ for $\epsilon=10^{-10}$, $T=25000$ and $\eta=0.7$.}
	\label{muoptfig}
	\centering
	\includegraphics[width=0.7\textwidth]{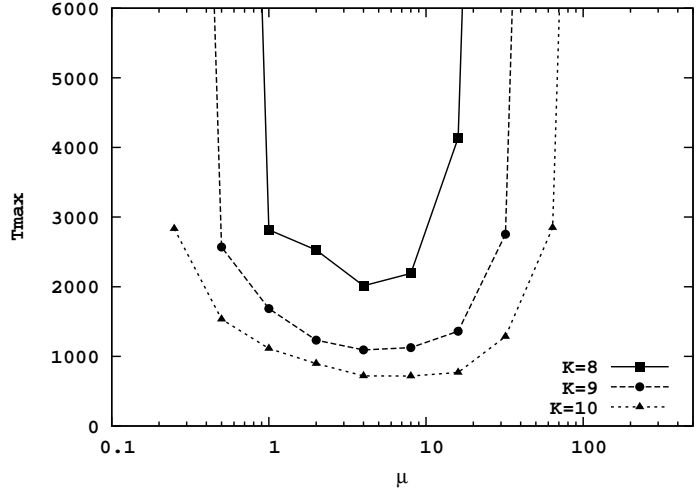}
\end{figure}

Figure \ref{muoptfig} plots the data in Table \ref{hmaxusmooth}.
From this figure it is clear that
the data
assimilation algorithm given by (\ref{alg2sys}) with $I_h$
replaced by $\tilde I_h$
works well when $\eta=0.7$ for a wide range
of values of the relaxation parameter $\mu$.

\section{Numerical Methods}

All fluid dynamics simulations presented in this paper were
performed
using a new parallel code written for the NVIDIA Compute Unified Device
Architecture in the CUDA C programming language \cite{nvidia2012}
which was developed on desktops at the University of Nevada Reno
and run on the Big Red II Cray XE6/XK7 supercomputer at
Indiana University.
Computations were made in the stream function formulation using
a fully-dealiased spectral Galerkin method.
Time steps were performed using a split-Euler method in which
the linear term was integrated exactly and the non-linear
terms were integrated using forward differences.

Specifically we set $\Delta\Psi=\curl u$ and compute the
reference solution using the stream function formulation
\begin{equation}\label{stream}
	\displaystyle {\partial \Delta \Psi\over \partial t}
	-\nu\Delta^2\Psi
	+\beta(\Psi) = \curl f,
\end{equation}
where
\begin{plain}$$\eqalign{
    \beta(\Psi)&= J(\Psi, \Delta \Psi)
		= \Psi_x \Delta \Psi_y -  \Psi_y \Delta \Psi_x\cr
    &=((\Psi_x)^2-(\Psi_y)^2)_{xy}
	    -(\Psi_x \Psi_y)_{xx}+(\Psi_x \Psi_y)_{yy}.\cr
}$$\end{plain}%
Similarly, set $\Delta\Phi=\curl v$ and compute
the approximating solution using
\begin{equation}\label{streamapp}
	\displaystyle {\partial \Delta \Phi\over \partial t}
	-\nu\Delta^2\Phi
	+\beta(\Phi) = \curl f -
	\mu (R_h(\Phi)-R_h(\Psi)),
\end{equation}
where $R_h(\Psi)= \curl P_\sigma I_h(\curl^{-1} \Delta \Psi).$
Note that $R_h\colon V^{3}\to V^{-1}$.

Following the 2/3 dealiasing rule applied to $512\times 512$ sized
discrete Fourier transforms we set $\K=\{-341,\ldots,341\}^2$ and
approximate
$$
    \Psi
     \approx
     \sum_{k\in\K}
        \hat{\Psi}_k \,e^{i k\cdot x},\quad
    \Phi
     \approx
     \sum_{k\in\K}
        \hat{\Phi}_k \,e^{i k\cdot x}
\words{and}
	\curl f = \sum_{k\in\K} \hat g_k e^{ik\cdot x}.
$$
Substituting these approximations into (\ref{stream})
and (\ref{streamapp})
and projecting onto the Fourier modes with wave numbers
$k\in\K$ yields the Galerkin truncations
$$
    - \displaystyle\frac{d\hat\Psi_k}{d t} k^2
    - \nu \hat{\Psi}_k k^4
    + \hat\beta(\Psi)_k = \hat g_k
$$
and
$$
    - \displaystyle\frac{d \hat\Phi_k}{d t} k^2
    - \nu \hat{\Phi}_k k^4
    + \hat\beta(\Phi)_k = \hat g_k -\mu \hat R_h(\Phi-\Psi)_k.
$$
The corresponding numerical scheme for the reference solution
is
\begin{plain}$$
	\displaystyle \hat \Psi_k(t+\Delta t)
		\approx e^{-\nu k^2 \Delta t}
		\Big\{\hat \Psi_k(t)
			+{\Delta t\over k^2} \hat\beta(\Psi(t))_k \Big\}
		-{1\over \nu k^4} \hat g_k \big( 1-e^{-\nu k^2 \Delta t}\big)
$$
and for the approximating solution is
$$\eqalign{
	\displaystyle \hat \Phi_k(t+\Delta t)
		&\approx e^{-\nu k^2 \Delta t}
		\Big\{\hat \Phi_k(t)
+{\Delta t\over k^2}\Big( \hat\beta(\Phi(t))_k
+\mu \hat R_h(\Phi-\Psi)_k\Big)\Big\}\cr
			&
\phantom{\,\approx e^{-\nu k^2 \Delta t}
        \Big\{\hat \Psi_k(t)
            +{\Delta t\over k^2} \hat\beta(\Psi(t))_k \Big\}\,}
		-{1\over \nu k^4} \hat g_k \big( 1-e^{-\nu k^2 \Delta t}\big).
}$$\end{plain}%

At the discrete level it is still the case that
only nodal-point observational measurements of the reference
solution are used to construct the approximating solution.
Moreover, since both solutions are integrated using the same
numerical methods, we may think of $\Psi$ as an unknown
reference solution that evolves according to a
known discrete dynamical system, and the solution represented
by $\Phi$ as an approximation generated by data assimilation
according to the exact same discrete dynamics.
Thus, even though our numerical schemes
are only first order in time, we consider our numerical
experiments to simulate data assimilation
in the absence of both measurement and model errors.

We take the time step to be $\Delta t=1/2048$ which is enough to
ensure that the CFL condition
$$
	{N \Delta t\over 2L} \sup_{x\in\Omega}\big\{
		|u_1(x)|+|u_2(x)|\big\} \le 0.1608\ll 1
$$
is satisfied for the reference solution over the entire run.
Note that as $\mu$ gets larger the data assimilation equations
(\ref{alg2sys}) become stiffer.  Therefore, the time step $\Delta t$
was also chosen small enough to ensure the stability of the
coupled numerical scheme for computing the approximating solution.

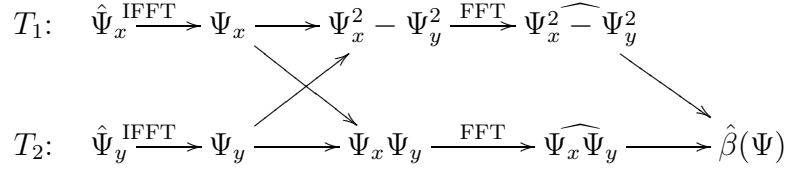
\begin{figure}[h!]
    \caption{Data Flow for Computing the Non-linear Term}
	\vskip-0.5cm
    \centering
    \begin{displaymath}
    \xymatrix{
		T_1\colon\quad
		\hat{\Psi}_{x} \ar[r]^{\rm IFFT} & \Psi_x \ar[dr]
        \ar[r] & \Psi_x^2 - \Psi_y^2 \ar[r]^{\rm FFT}
        & \widehat{\Psi_x^2 - \Psi_y^2} \ar[rd]  \\
		T_2\colon\quad
        \hat{\Psi}_{y} \ar[r]^{\rm IFFT} & \Psi_y
        \ar[r] \ar[ru] & \Psi_x \Psi_y \ar[r]^{\rm FFT}
        & \widehat{\Psi_x \Psi_y} \ar[r]
        & \hat\beta(\Psi)
    }
    \end{displaymath}
	\vskip-0.5cm
    \label{2d_fft_num}
\end{figure}

Our numerical software has been optimized so that it runs entirely
on the CUDA hardware with zero
memory copies and four fast Fourier transforms per time step.
Note that four transforms is the minimal number
for the two-dimensional Navier--Stokes equations,
see Basdevant \cite{basdevant1983} for further remarks and
analogous optimizations when computing the
three-dimensional Navier--Stokes equations.
As device memory on the CUDA hardware is relatively scarce
we also minimized our storage requirements.
Storage requirements consist of four $n\times n$
double-precision scalar arrays: one for $\Psi$,
another for $\Phi$ and two temporary arrays
$T_1$ and $T_2$.
Figure~\ref{2d_fft_num} shows the data flow
when computing the non-linear term.
The first line represent the contents of $T_1$,
the second represents $T_2$ and the arrows represent
computational kernels.

When using
$512\times 512$ FFTs, our code achieves approximately 1062 time steps
per second.
In particular, the computational speed running on an NVIDIA Tesla K20 GPU
was found to be roughly 37 times faster than equivalent CPU code
running on single AMD Opteron 6212 core and 11.5 times faster
when compared to running on 32 CPU cores.
Correctness of operation was verified using the
Navier--Stokes solver described in~\cite{olson2003} and~\cite{olson2008}.

\section{Conclusion}

As is consistent with the analytical bound (\ref{alg2thm})
and related discussion in \cite{azouani2014}, the numerical
results given in Table \ref{hmaxu} and Figure \ref{muoptfig}
show that the approximating solution does not converge to the
reference solution when $\mu$ is either too small or too big.
At the same time, provided the resolution $h$ is fine enough
and $\eta\approx 0.7$, there
is a wide range of good values for $\mu$ when using
the smoothed interpolant observable $\tilde I_h$.
In particular, when
$h=L/10\approx 0.6283$, the data assimilation algorithm
(\ref{alg2}) performs similarly for values of $\mu$
between $0.5$ and $32$.
Note, however, that smaller values of $\mu$ are computationally
preferable because of stiffness considerations.

We remark that our numerical experiments have been conducted
using exact error-free measurements and exact model dynamics
and that in the presence of measurement and model errors we don't
expect a similarly wide range of good values for $\mu$.
In fact, preliminary computations show when is noise added to
the system
that there exists a unique optimal value for $\mu$
reflecting the tradeoff between measurement and model errors.
Theoretically, if the dynamics represented by ${\cal F}(u)$ in
(\ref{refsol}) are linear,
then $\mu$ can be seen as the parameter of a linear Kalman filter,
see for example Majda and Harlim \cite{majda2012},
and there exists an analytically derived optimal value
for $\mu$ which represents the tradeoff between measurement
and model errors.
In the fully non-linear case studied here, the fact that
there is a wide range of good constant values for $\mu$
in the absence of measurement and model errors
suggests, provided the
resolution $h$ is fine enough, that $\mu$ can be further
optimized as if the model was linear.

We now compare the coarsest resolution $h=L/K$ that
works for the data assimilation experiments presented here
when $\eta=0.7$ to the number of numerically determining modes
found in
\cite{olson2008}.
Under the same physical parameters and forcing, the
minimum number
of Fourier modes needed by (\ref{alg1})
was
$$
	n_c={\rm card} ({\cal D}_5)=80$$
where
$${\cal D}_{R}
	=\{\,e^{i k\cdot x}: 0<k_1^2+k_2^2\le R^2\,\}.$$
In this paper we show the minimum of nodal
measurements needed are
$$
	N=K^2=64
$$
which by the Nyquist--Shannon sampling theorem may be
represented by the Fourier modes
$${\cal N}_K=
	\{\,e^{i k\cdot x}: 0<
		 \max(|k_1|,|k_2|)\le K/2\,\}.$$

To compare these two results we note that
${\cal D}_R$ represents a circle in Fourier space
while ${\cal N}_K$ represents a square.
Let $R_{\rm min}=5$ and $K_{\rm min}=8$.
If the resolution requirements of algorithm (\ref{alg1}) are
comparable to (\ref{alg2}), we would expect that
${\cal N}_{K_{\rm min}}\subseteq {\cal D}_R$
would imply $R\ge R_{\rm min}$ and that
${\cal D}_{R_{\rm min}}\subseteq {\cal N}_K$
would imply $K\ge K_{\rm min}$.
This is supported by our results.
If
${\cal N}_{8}\subseteq {\cal D}_R$ then
$$R\ge {8\sqrt 2\over 2}\approx 5.65 \ge 5=R_{\rm min}.$$
Similarly, if
${\cal D}_{5}\subseteq {\cal N}_K$ then
$$K\ge 2\cdot 5=10\ge 8=K_{\rm min}.$$

Thus, even though our nodal observations possess the problems
of aliasing and high-frequency spill over, these problems
can be mediated with appropriate smoothing.  The resulting
resolution $K$
needed for the approximating solution to converge to the
exact solution is then about the same as suggested by the
number of numerically determining modes.

\section*{Acknowledgements}

The authors would like to thank Professor Michael Jolly for his
help with the Big Red II supercomputer at Indiana and his current
collaboration on research treating measurement and model error.
The work of E.O. was supported in part by EPSRC grant EP/G007470/1,
by sabbatical leave from the University of Nevada Reno and by
NSF grant DMS-1418928.
The work of E.S.T. was supported in part by  a grant of the ONR, and the NSF
grants DMS-1109640 and DMS-1109645.

\vfill\eject

\end{document}